\documentclass[a4paper]{article}
\usepackage{spconf,amsmath,graphicx,cite}

\def\x{{\mathbf x}}
\def\t{{\mathbf t}}

\def\dkt{{\textrm{d}^k {\mathbf t}}}

\def\R{{\mathbf{R}}}
\def\IMF{{\textrm{IMF}}}

\title{Multidimensional Iterative Filtering method for the decomposition of high--dimensional non--stationary signals}
\twoauthors
  {Antonio Cicone\sthanks{This work was supported by the Italian ``National Group for Scientific Computation'' (GNCS-INDAM) - Programma Giovani Ricercatori 2014}}
	 {DISIM\\
      Universit\'a degli Studi dell'Aquila\\
      via Vetoio 1, 67100, L'Aquila, Italy}
  {Haomin Zhou\sthanks{This work was supported by NSF Faculty Early Career Development
(CAREER) Award DMS--0645266, DMS--1042998, DMS--1419027, and ONR Award
N000141310408.}}
	{School of Mathematics\\
	Georgia Institute of Technology\\
	686 Cherry St., Atlanta, GA 30332, USA}

\begin{document}

\maketitle
\begin{abstract}
Iterative Filtering (IF) is an alternative technique to the Empirical
Mode Decomposition (EMD) algorithm for the decomposition
of non--stationary and non--linear signals. Recently in \cite{cicone2014adaptive} IF has been
proved to be convergent for any $L^2$ signal and its stability
has been also showed through examples. Furthermore in \cite{cicone2014adaptive} the so called
Fokker--Planck (FP) filters have been introduced. They are
smooth at every point and have compact supports. Based on
those results, in this paper we introduce the Multidimensional
Iterative Filtering (MIF) technique for the decomposition and
time--frequency analysis of non--stationary high--dimensional signals.
And we present the extension of FP filters to higher dimensions.
We illustrate the promising performance of MIF
algorithm, equipped with high--dimensional FP filters, when
applied to the decomposition of 2D signals.
\end{abstract}
\section{Introduction}
\label{sec:intro}
Given a non--stationary signal a hard problem is how to perform a time--frequency analysis in order to unravel its hidden features. The problem becomes even harder if we want to handle a signal that have dimension higher than one. Such kind of problems are ubiquitous in real life and their solutions can help shedding lights in many research fields like, for instance, in the non--destructive detection of structural damages in buildings and machineries, in the quantitative canvas weave analysis in the art investigations of paintings, in the identification of harmful airborne chemical particles by means of hyperspectral images analysis, in the image enhancement for medical applications, in the extraction of information from atomic crystal images, et cetera.

Time--frequency analysis have been substantially studied in the past \cite{cohen1995time} and, traditionally, linear techniques like Fourier spectral analysis or wavelet transforms have been commonly used to decompose signals, even non--stationary ones, into simple and stationary components and then perform on them a frequency analysis. These approaches, even if easy to implement, have limitations. In particular, they work well when the given signal is periodic and stationary, whereas they cannot handle properly non--stationary signals. Furthermore all these techniques use predetermined bases and they are not designed to be data--adaptive. Hence there is the need to develop new methods designed to handle specifically non--linear and non--stationary signals.

The research on these new techniques started in 1998 when Huang and his research group developed and released the very first method of this kind, the so called  Empirical Mode Decomposition (EMD) algorithm \cite{huang1998empirical}. The idea behind this method is to generalize the approach of classical linear techniques: decompose a signal into simple components and then perform a time--frequency analysis on each component separately. For this reason Huang and his group defined the so called Intrinsic Mode Functions (IMFs) which are simple functions with the following properties: their number of extrema is either equal to the number of zero crossings or they differ at most by one, and at any point their moving average is zero. Furthermore they developed an iterative technique, called the sifting process, to decompose a signal into such Intrinsic Mode Functions with the final goal of computing the instantaneous frequency \cite{boashash1992estimating,cohen1995time} of these simple components.

The sifting process of the EMD is structured in the following way. Let $L$ be an operator capturing the moving average of a signal $f$,
and let $S(f)(t) = f(t) - L(f)(t)$ be its fluctuation part. Then the first IMF is given by $\IMF_1 (t)= \lim_{n\rightarrow\infty}S_{1,n}(f)(t)$. The limit is such that the moving average of $\IMF_1 (t)$ is the zero function. Assuming $k\geq1$ IMFs have been computed, then the subsequent IMFs are produced applying the aforementioned procedure to the residual $r(t)=f(t)-\sum_{j=1}^k \IMF_j(t)$. The method stops when the residual $r$ becomes a trend signal. So in the end the given dataset is decomposed as $f(t) = \sum_{j=1}^m \IMF_j(t) + r(t)$. In this algorithm the moving average $L(f)(t)$ is given by the mean function of an upper and a lower envelope, where upper and lower envelopes are cubic splines connecting local maxima and local minima of $f(t)$, respectively. This method has received a lot of attention in the last decade. In fact,  many researchers have applied it to solve several open problems in diverse research fields.

We point out that there are known examples showing that this technique can be unstable to small perturbations of the initial data. To overcome this issue Huang et al. developed the Ensemble Empirical Mode Decomposition (EEMD) \cite{wu2009ensemble} where the IMFs are taken as the mean of many different trials produced using the EMD algorithm. Since cubic splines are used repeatedly in the iterations of both EEMD and EMD, a rigorous proof of the convergence of these methods is still missing.

Inspired by the seminal work by Huang et al., in the last years many research groups started working to the development of alternative techniques to the EMD algorithm. For instance Daubechies and her research group devised the Synchrosqueezed wavelet transform \cite{daubechies2011synchrosqueezed}, Gilles created the Empirical wavelet transform \cite{gilles2013empirical}, Hou et al., using the multicomponent amplitude modulation and frequency modulation (AM--FM) representation \cite{wei1998instantaneous}, developed the Sparse time--frequency representation method \cite{hou2011adaptive}, and many others \cite{pustelnik2012multicomponent,selesnick2011resonance, dragomiretskiy2014variational,wu2011one}. All these alternative methods make use of optimization techniques for the decomposition of a given non--stationary and non--linear signal. These approaches require the a priori selection of a suitable basis for the decomposition.

In 2009 Zhou et al. developed another alternative algorithm called Iterative Filtering \cite{lin2009iterative}. This method, as the original EMD, has an iterative structure and, unlike all the techniques based on optimization, it does not require any initial assumption on the signal. Therefore it is able to produce decompositions that are completely data driven. This method has the same structure of the sifting process, but now the moving average $L(f)(t)$ is derived by convolution of the given signal $f(t)$ with filters like, for instance, a double average filter. This new way of computing moving averages allows to perform a rigorous analysis of this technique.
In particular, under mild sufficient conditions on the filters used in the convolutions, the analytic convergence of Iterative Filtering applied to generic $L^2$ signals is ensured and the components produced in the decomposition of a signal can be derived from an explicit analytic formula \cite{cicone2014adaptive}.

Some of the aforementioned methods, initially developed to decompose 1D datasets, have been already generalized to handle two and higher dimensional signals, like for instance the multidimensional ensemble empirical mode decomposition method \cite{wu2009multi} or the Synchrosqueezed wave packet and curvelet transform \cite{yang2014synchrosqueezed}. However this is not the case of Iterative Filtering. Therefore in this paper we address the problem of extending Iterative Filtering and the so called Fokker--Planck (FP) filters \cite{cicone2014adaptive} to higher dimensions. We introduce, in particular, the Multidimensional Iterative Filtering algorithm and the Generalized Fokker--Planck (GFP) filters. We show also the performance of this method, equipped with GFP filters, on artificial and real life signals.

\section{Multidimensional Iterative Filtering algorithm and generalized Fokker--Planck filters}
\label{sec:MIF}

The \emph{Multidimensional Iterative Filtering} (MIF) technique consists of an Inner and an Outer Loop. In the inner loop the method computes an IMF of a given $k$--dimensional signal $f$ as the limit of the sequence generated by subtracting from the signal its moving average iteratively. In the outer loop we simply update the signal by removing from it the previously computed IMFs. The outer loop is iterated until the remainder becomes a trend signal. Whereas the inner loop should be theoretically iterated until the moving average becomes a zero function. However, in the numerical implementation of the inner loop, we use a stopping criterion to discontinue the iterations. The pseudocode of this algorithm is given in Table~\ref{tab:MIF}, where $f$ is the $k$--dimensional signal we want to decompose and $w\in\R^k$ is a filter function with finite support $\Omega\subset\R^k$.

\begin{table}
\small
\begin{center}
\renewcommand{\arraystretch}{1.2}
\begin{tabular}{l}
\hline
\textbf{ MIF Algorithm } IMFs = MIF$(f)$ \\
\hline
IMFs = $\left\{\right\}$ \\
\textbf{while} the average number of extrema of $f$ $\geq 2$ \textbf{do} \\
$\quad$ compute the filter support $\Omega$ for $f$\\
$\quad$ $f_1 = f$\\
$\quad$ \textbf{while} the stopping criterion is not satisfied \textbf{do}\\
$\quad\quad$ $f_{n+1}(\x) = f_{n}(\x)-\int_{\Omega} f_n(\x+\t)w(\t)\dkt$\\
$\quad\quad$ $n = n+1$\\
$\quad$ \textbf{end while}\\
$\quad$ IMFs = IMFs$\,\cup\,  \{ f_{n}\}$\\
$\quad$ $f=f-f_{n}$\\
\textbf{end while}\\
IMFs = IMFs$\,\cup\,  \{ f\}$\\
\hline
\end{tabular}
\end{center}
\caption{Multidimensional Iterative Filtering pseudocode}\label{tab:MIF}
\end{table}

As a filter $w$ we need to have, first of all, nonnegative $L^2$ functions which are finitely supported on $\Omega\subset\R^k$ with $\int_{\Omega} w(\t)\dkt = 1$. However not every $w$ function which fulfils the previous properties is going to ensure the convergence of this technique. A well chosen class of filters has to be adopted. From what is known for the 1D Iterative Filtering algorithm \cite{cicone2014adaptive}, we can conjecture that in higher dimensions a class of filters which guarantees the MIF convergence is given by a proper extension of the 1D Fokker--Planck (FP) filters. Such 1D filters, derived from the solution of Fokker--Planck partial differential equations \cite{cicone2014adaptive}, have the nice property of being extremely smooth since they are infinitely differentiable at any point of their domain. To produce axial symmetric Fokker--Planck filters of dimension higher than one there is no need to numerically solve high--dimensional partial differential equations. We can make use of highly accurate numerical solutions in 1D to generate, first, numerical 2D Fokker--Planck filters and then higher dimensional ones. This can be achieved by simply scaling, resampling, and then rotating around one axis a 1D filter many times to produce a 2D version of it. Higher dimensional filters can be produced simply iterating this procedure. We call them \emph{Generalized Fokker--Planck} (GFP) filters. We observe that this approach allows to produce also GFP filters with an ellipsoidal support $\Omega$.

Once a filter shape has been selected it can be used for every dataset. The problem is, in order to extract meaningful IMFs, how to select a proper support $\Omega$ for the filter simply based on the signal we want to decompose. Following what proposed in \cite{lin2009iterative} for the 1D case, we compute for each dimension of the signal the average 1D length of the support. Then we use this information to either find the radius of a spherical support or to identify the radii of an ellipsoidal $\Omega\subset\R^k$.

About the stopping criterion, there are many possibile options. One way of doing it, as mentioned in \cite{cicone2014adaptive,huang1998empirical} for the 1D case, is to consider the relative change in $f_{n}$ and discontinue the inner loop as soon as a prefixed threshold is reached.

Finally, we observe that this extension to higher dimensions can be applied in the same way also to the so called Adaptive Local Iterative Filtering (ALIF) method \cite{cicone2014adaptive}, which is a 1D generalization of the Iterative Filtering algorithm.

\section{Numerical examples}
\label{sec:ex}

In the following we show the performance of the Multidimensional Iterative Filtering (MIF) algorithm, equipped with Generalized Fokker--Planck (GFP) filters, applied to both artificial and real life signals. For simplicity from now on we consider only bidimensional signals, however we point out that the MIF method can handle datasets of any dimension.

\subsection{Example 1}

\begin{figure}
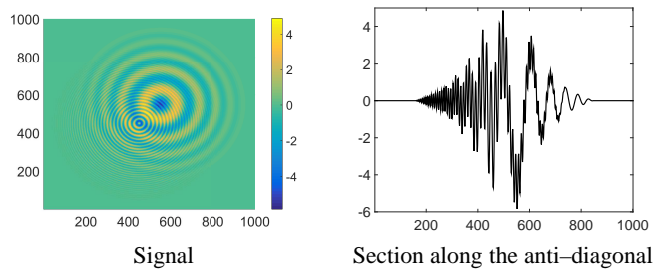

\centering
\begin{minipage}[b]{.48\linewidth}
  \centering
  \centerline{\includegraphics[width=\linewidth]{Ex1_Signal.eps}}
  \small\centerline{Signal}\medskip
\end{minipage}
\hfill
\begin{minipage}[b]{0.45\linewidth}
  \centering
  \centerline{\includegraphics[width=\linewidth]{Ex1_Signal_sec.eps}}
  \small\centerline{Section along the anti--diagonal}\medskip
\end{minipage}
\vspace*{-0.3cm}
\caption{Artificial bidimensional dataset containing a mixture of two non--stationary sinusoidal signals.}
\label{fig:Ex1Sig}
\end{figure}

We start considering the artificial signal showed in Figure~\ref{fig:Ex1Sig} containing the mixture of two non--stationary sinusoidal damped signals.

\begin{figure}
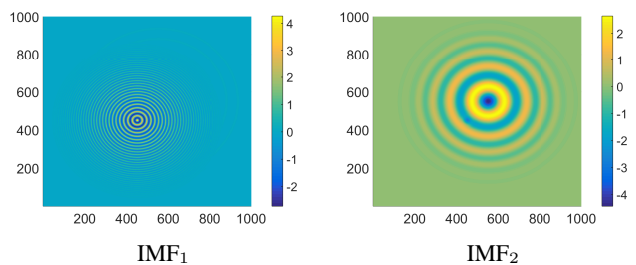

\centering
\begin{minipage}[b]{.48\linewidth}
  \centering
  \centerline{\includegraphics[width=\linewidth]{Ex1_IMF1.eps}}
  \small\centerline{$\IMF_1$}\medskip
\end{minipage}
\hfill
\begin{minipage}[b]{0.48\linewidth}
  \centering
  \centerline{\includegraphics[width=\linewidth]{Ex1_IMF2.eps}}
  \small\centerline{$\IMF_2$}\medskip
\end{minipage}
\vspace*{-0.3cm}
\caption{IMFs produced using MIF.}
\label{fig:Ex1IMFs}
\end{figure}

If we apply the MIF algorithm to this signal we obtain the two IMFs plotted in Figure~\ref{fig:Ex1IMFs}.

\begin{figure}
\centering
\begin{minipage}[b]{.48\linewidth}
  \centering
  \centerline{\includegraphics[width=\linewidth]{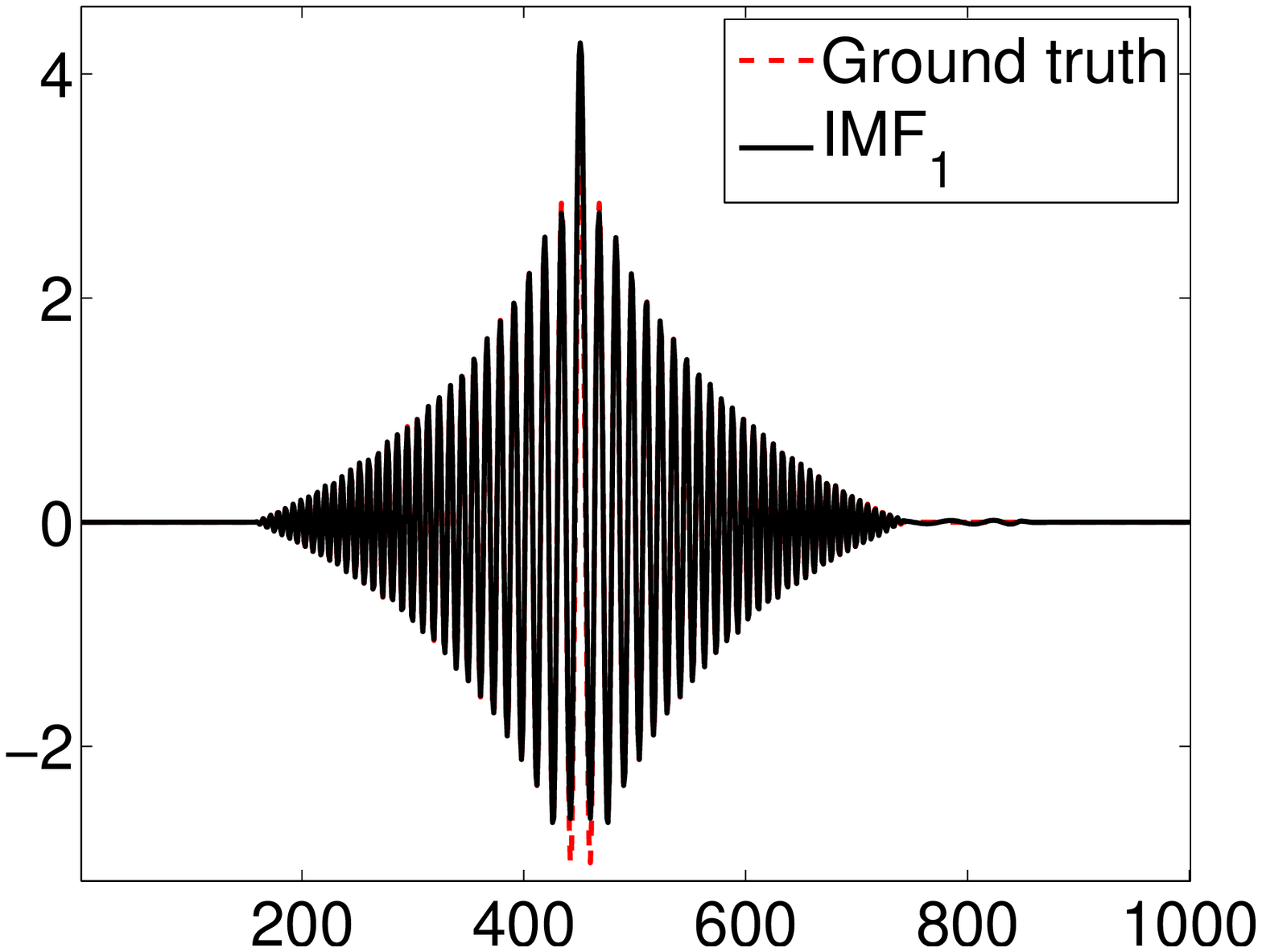}}
  \small\centerline{$\IMF_1$ vs ground truth}\medskip
\end{minipage}
\hfill
\begin{minipage}[b]{0.48\linewidth}
  \centering
  \centerline{\includegraphics[width=\linewidth]{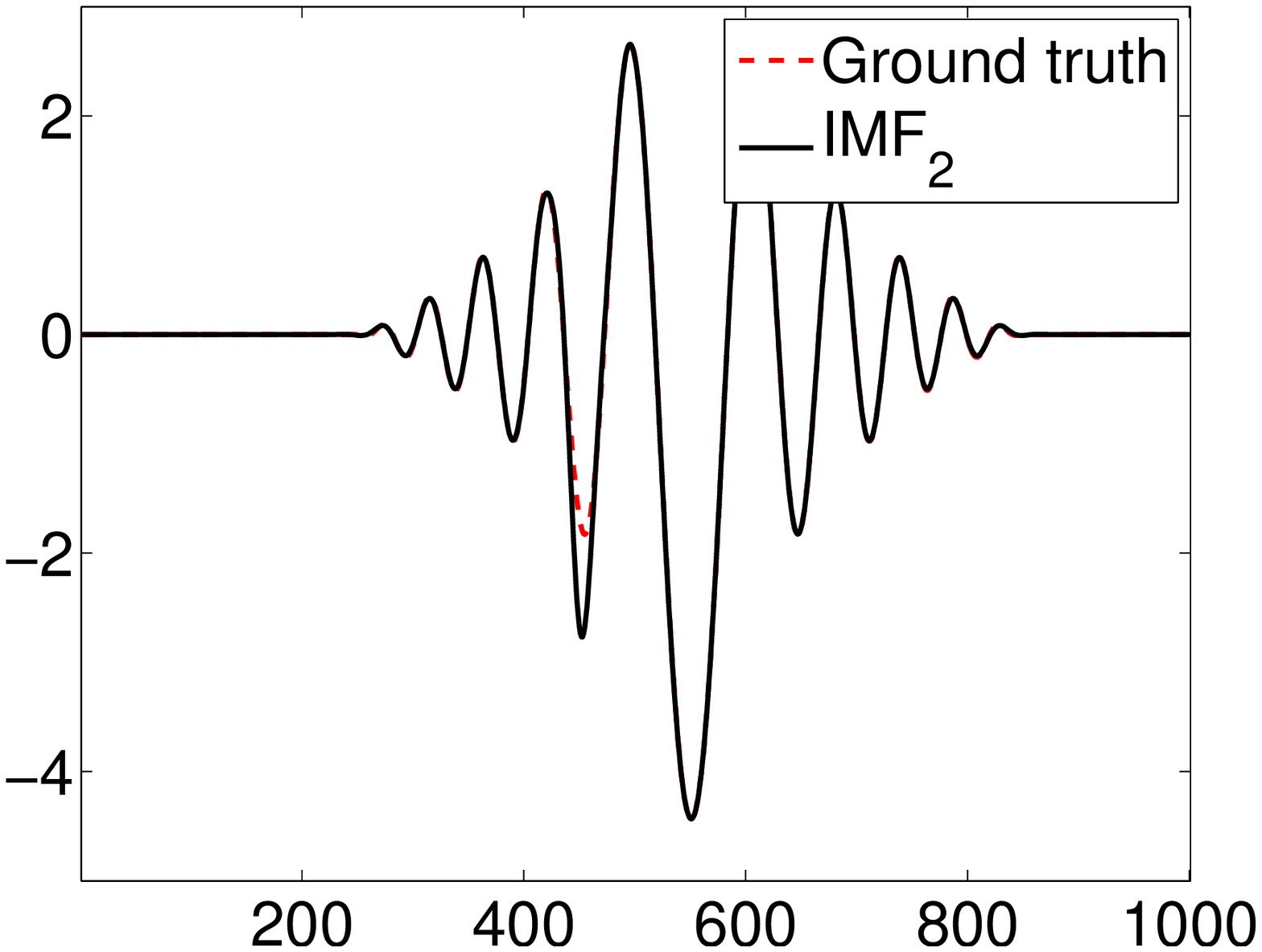}}
  \small\centerline{$\IMF_2$ vs ground truth}\medskip
\end{minipage}
\vspace*{-0.3cm}
\caption{Sections along the anti--diagonal of the IMFs.}
\label{fig:Ex1_sec}
\end{figure}

To have a better understanding of the accuracy of the decomposition we can plot sections of these IMFs taken along anti--diagonals. This approach allows to easily compare the IMFs with the ground truth, as showed in Figure~\ref{fig:Ex1_sec}.

\begin{figure}
\centering
\begin{minipage}[b]{.48\linewidth}
  \centering
  \centerline{\includegraphics[width=\linewidth]{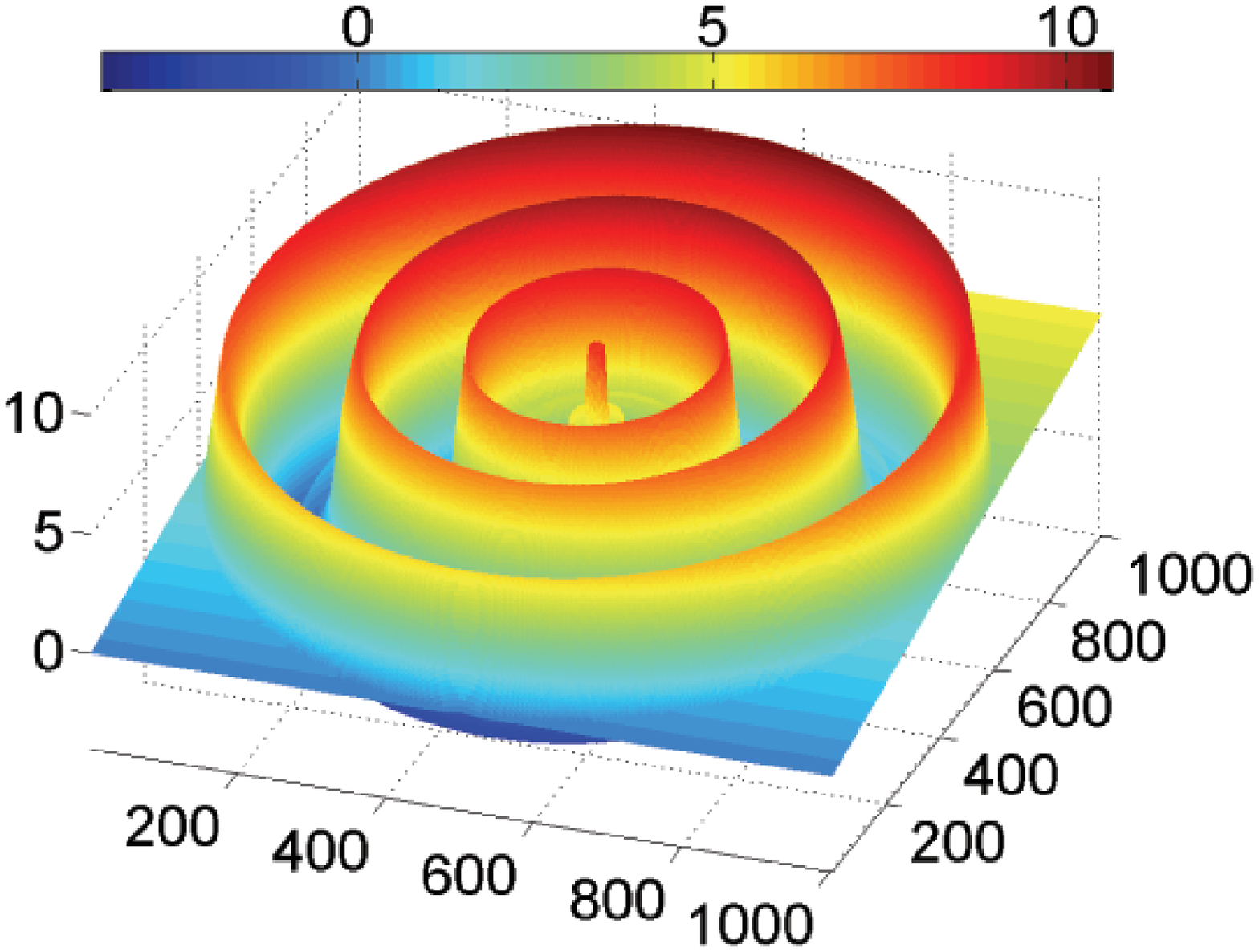}}
\end{minipage}
\hfill
\begin{minipage}[b]{0.45\linewidth}
  \centering
  \centerline{\includegraphics[width=\linewidth]{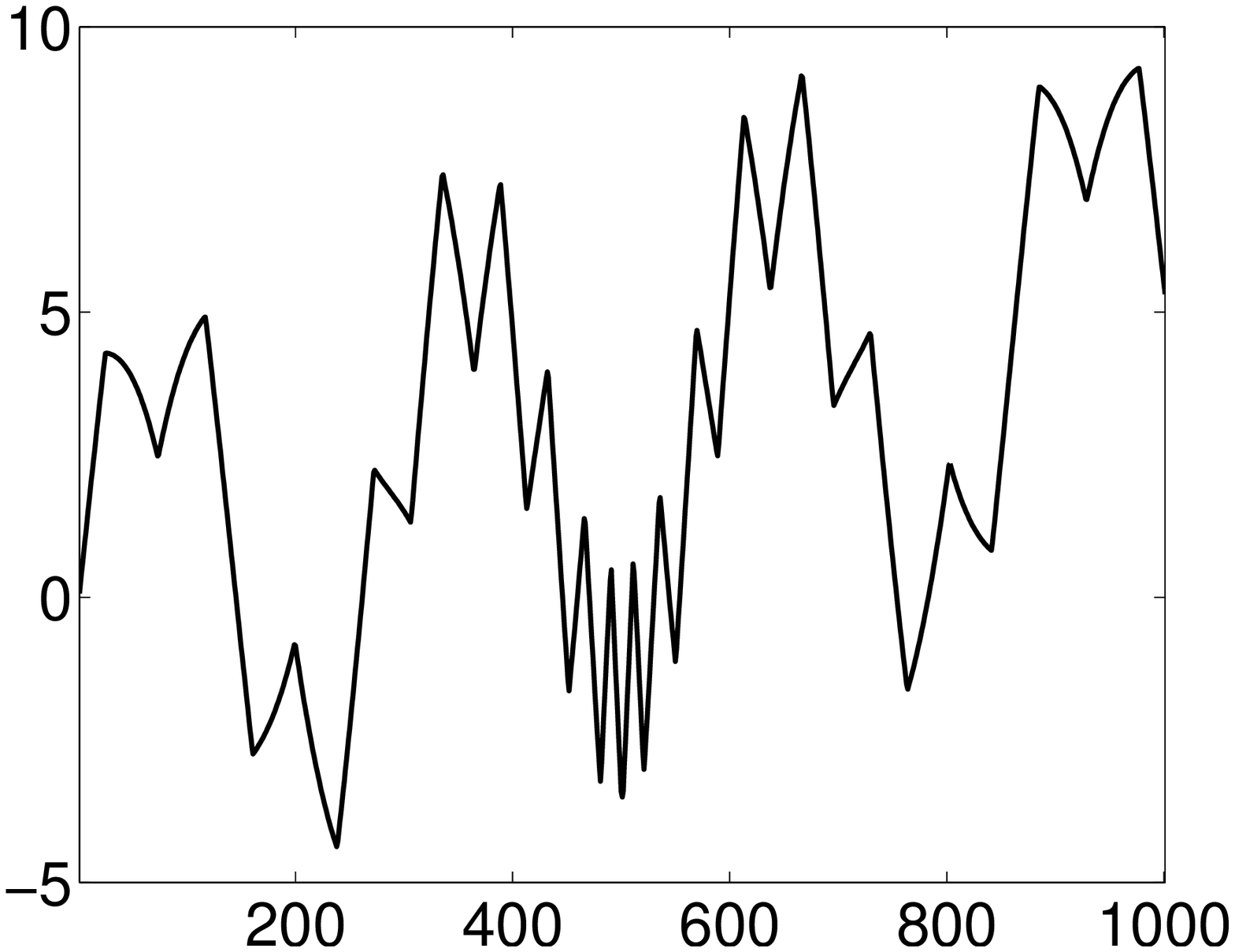}}
\end{minipage}
\vspace*{-0.3cm}
\caption{A non--stationary and non--smooth signal and its middle vertical section.}\label{fig:Ex2Sig}
\end{figure}

\subsection{Example 2}

\begin{figure}
\centering
\begin{minipage}[b]{.48\linewidth}
  \centering
  \centerline{\includegraphics[width=\linewidth]{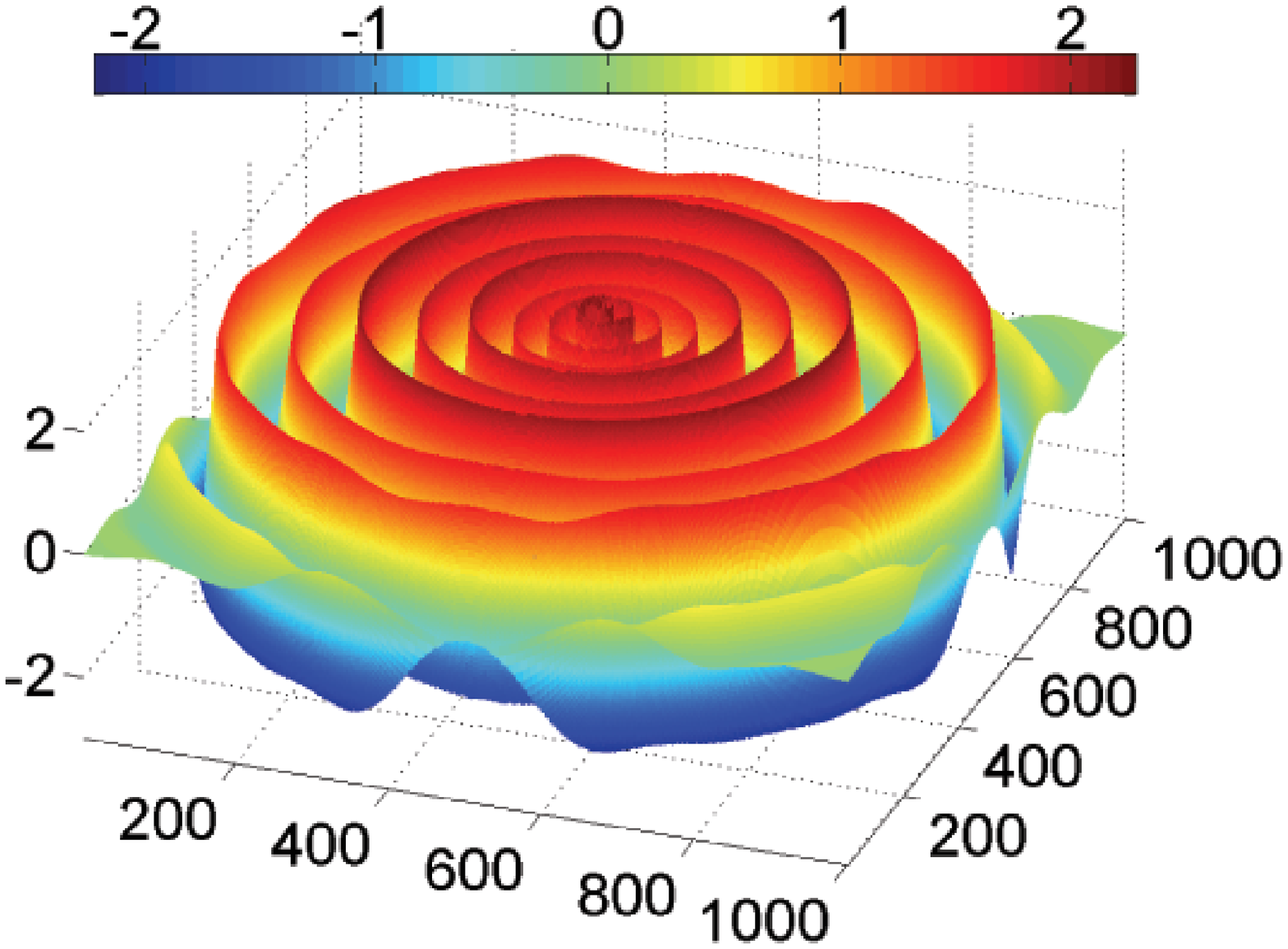}}
\end{minipage}
\hfill
\begin{minipage}[b]{0.48\linewidth}
  \centering
  \centerline{\includegraphics[width=\linewidth]{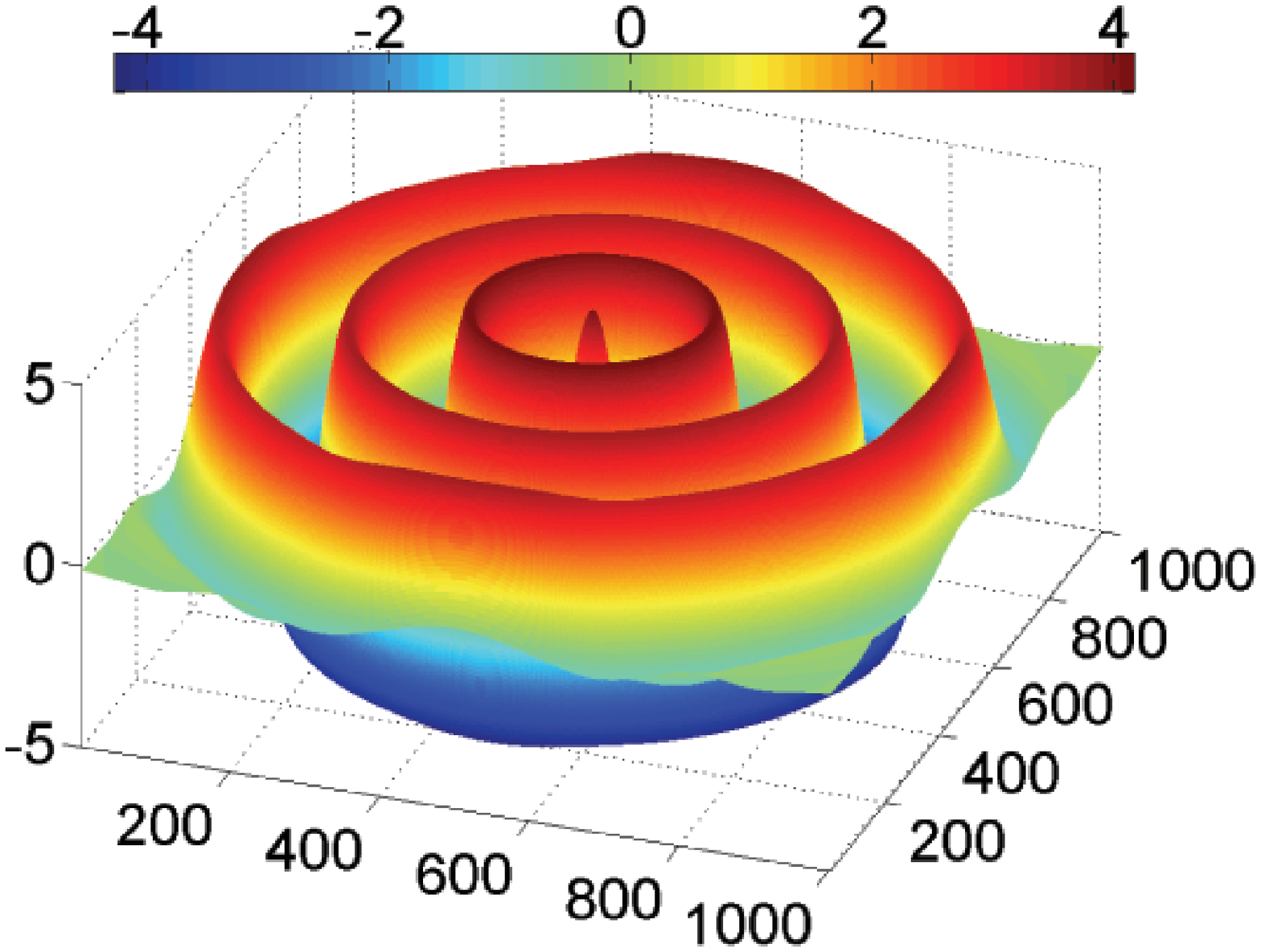}}
\end{minipage}
\vspace*{-0.3cm}
\caption{First (left) and second (right) IMF produced by MIF.}\label{fig:Ex2IMFs}
\end{figure}

\begin{figure}
\centering
\begin{minipage}[b]{.48\linewidth}
  \centering
  \centerline{\includegraphics[width=\linewidth]{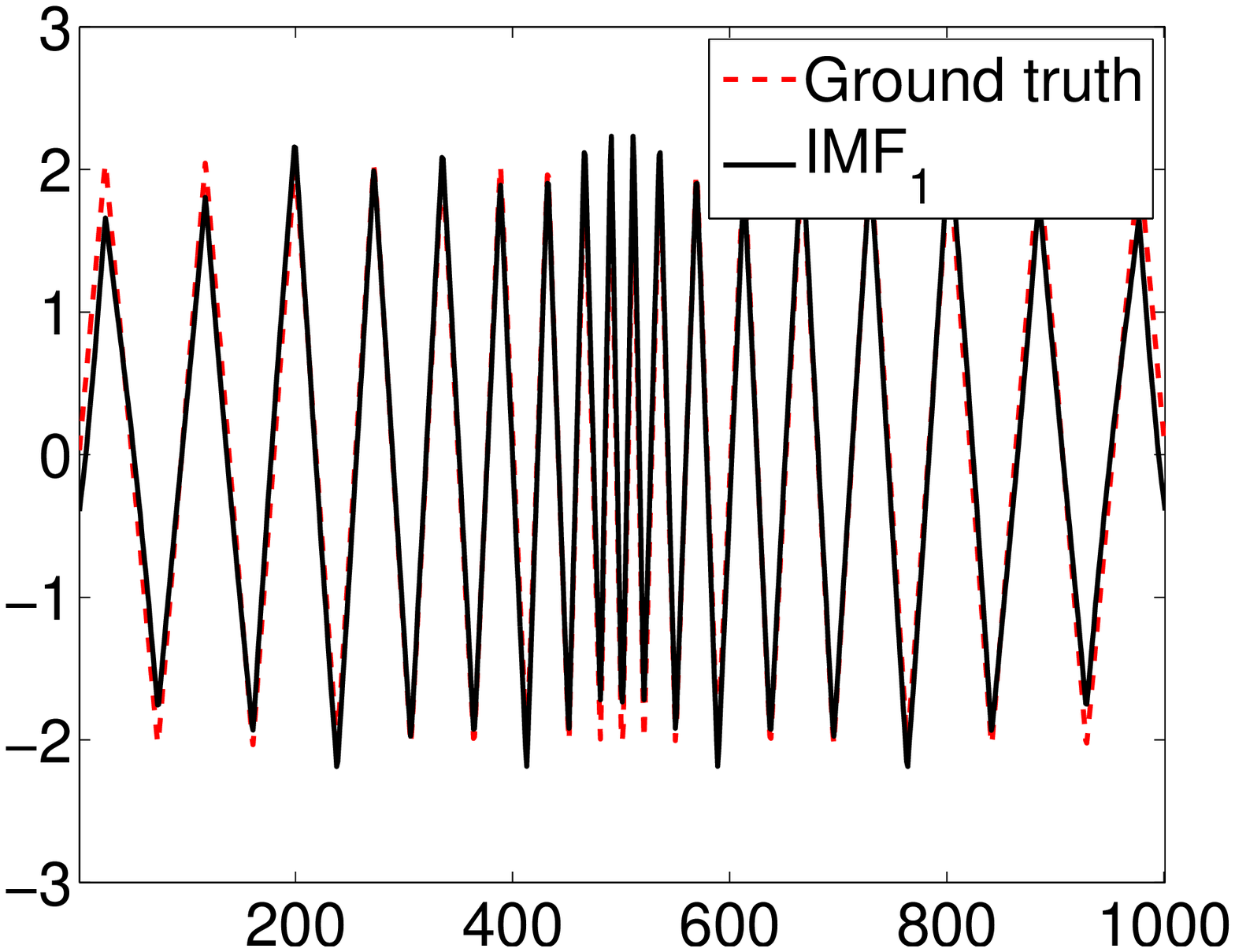}}
  \small\centerline{$\IMF_1$}\medskip
\end{minipage}
\hfill
\begin{minipage}[b]{0.48\linewidth}
  \centering
  \centerline{\includegraphics[width=\linewidth]{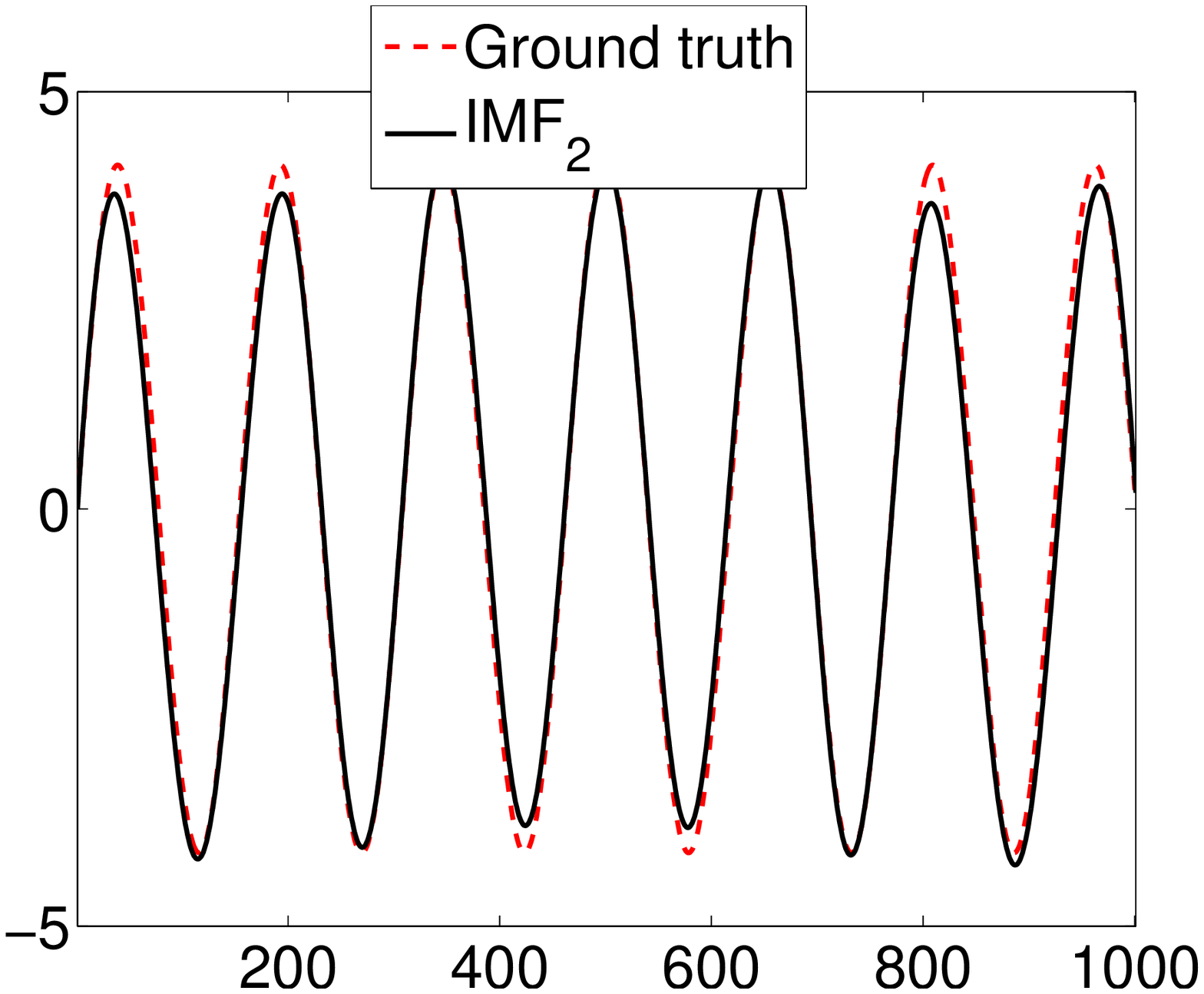}}
  \small\centerline{$\IMF_2$}\medskip
\end{minipage}
\vspace*{-0.3cm}
\caption{Middle vertical sections of the IMFs and the corresponding ground truth.}
\label{fig:Ex2_sec}
\end{figure}

\begin{figure}
\centering
\begin{minipage}[b]{.48\linewidth}
  \centering
  \centerline{\includegraphics[width=\linewidth]{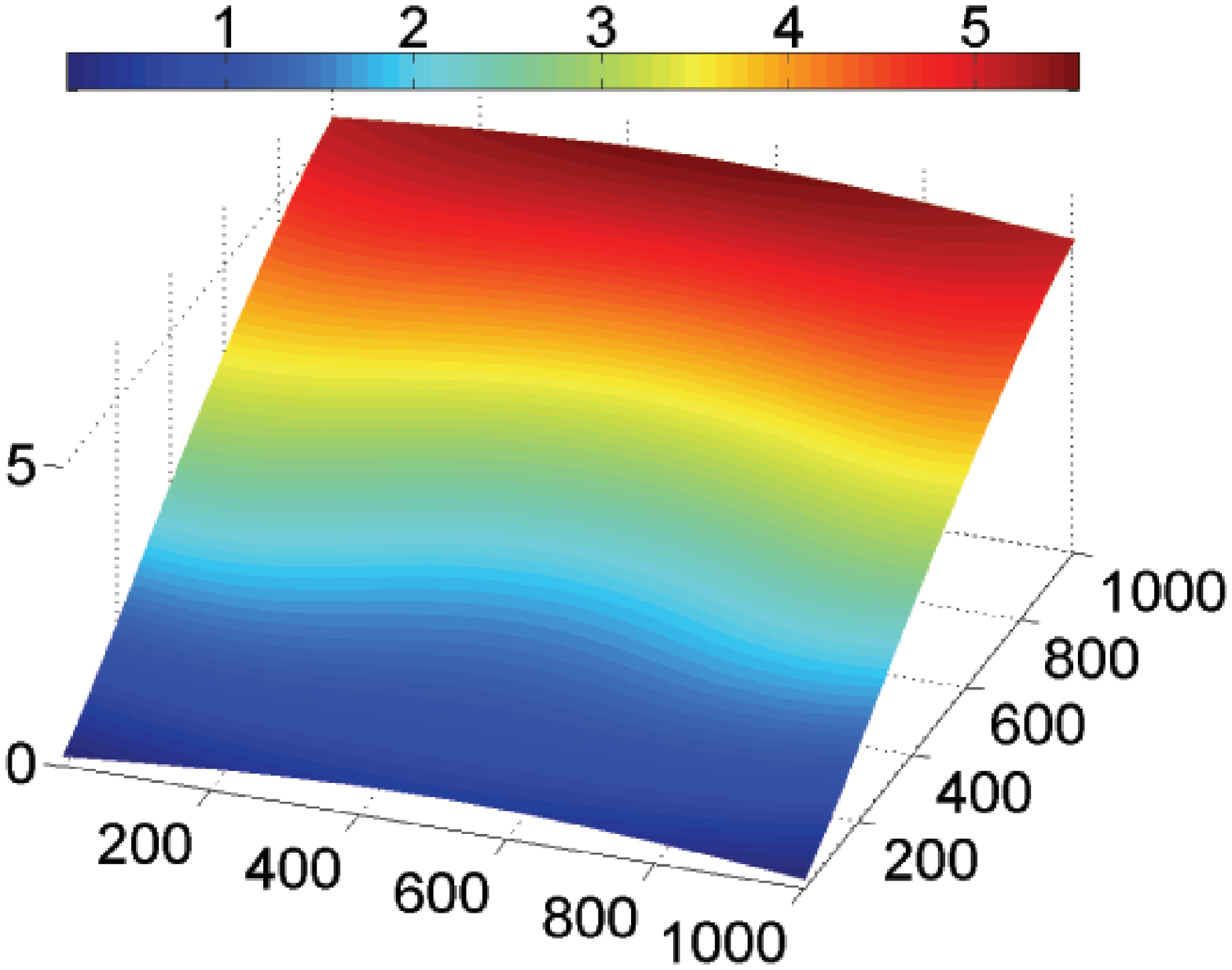}}
\end{minipage}
\hfill
\begin{minipage}[b]{0.48\linewidth}
  \centering
  \centerline{\includegraphics[width=\linewidth]{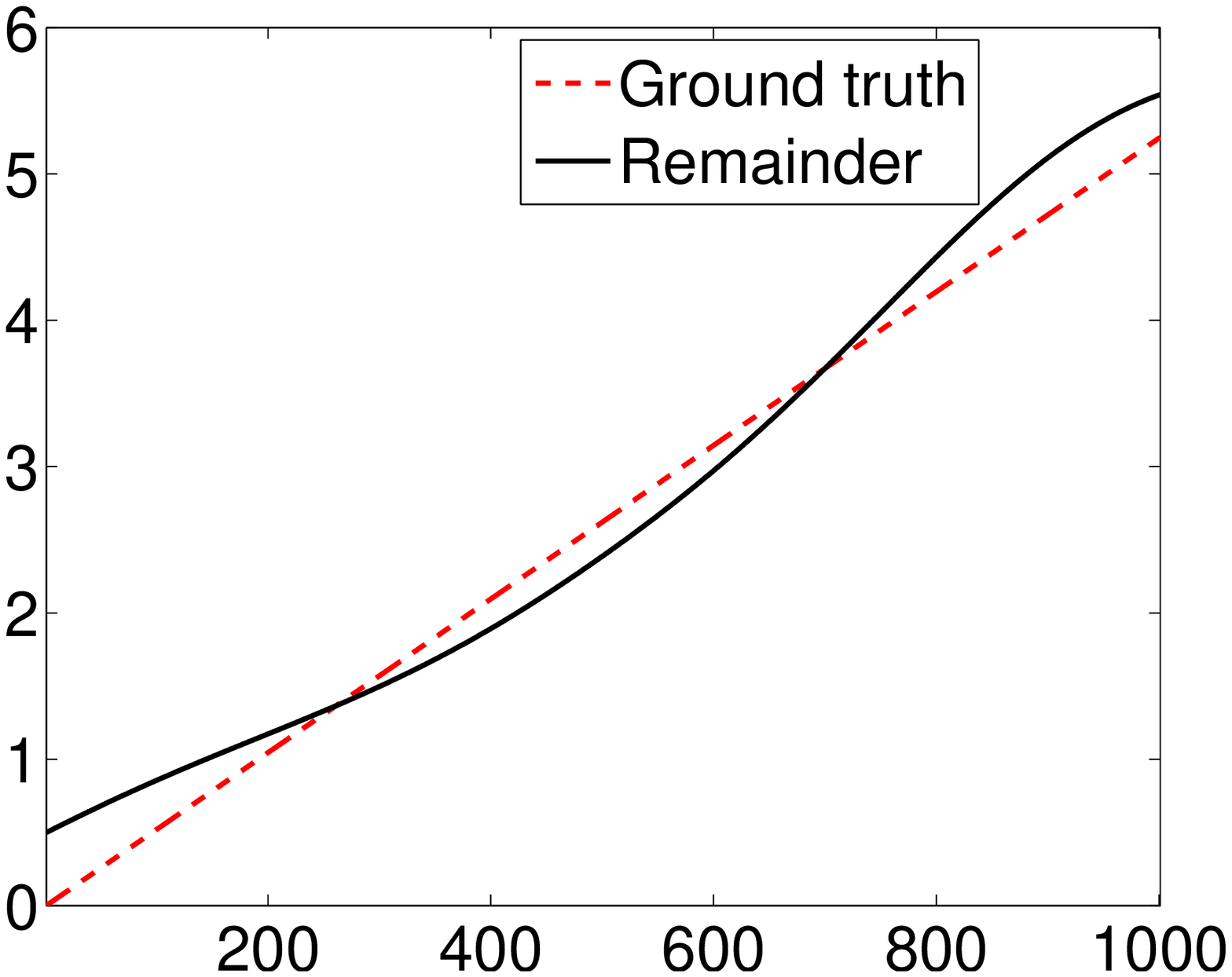}}
\end{minipage}
\vspace*{-0.3cm}
\caption{The Remainder, its middle vertical section and the corresponding ground truth.}
\label{fig:Ex2_Rem}
\end{figure}

For the second example we consider the signal of Figure~\ref{fig:Ex2Sig} which is the mixture of a smooth sinusoidal signal with a non--smooth and non--stationary one. If we run the method, using a 2D generalized Fokker--Planck filter, we can produce a first IMF which is a non--smooth and non--stationary signal and a second one which is smooth and stationary, as showed in Figure~\ref{fig:Ex2IMFs}. To have a better understanding of the kind of signal we are dealing with and to appreciate the performance of the decomposition algorithm, we plot middle vertical slices of the original signal, IMFs, remainder and ground truth in Figures~\ref{fig:Ex2Sig}, \ref{fig:Ex2_sec} and \ref{fig:Ex2_Rem}. This example shows the ability of this method to separate components completely different in nature, like smooth and non--smooth IMFs, using a fixed filter function. It would be hard to produce similar results with any optimization technique.

\subsection{Example 3}

\begin{figure}
\centering
\begin{minipage}[b]{0.8\linewidth}
  \centering
  \centerline{\includegraphics[width=\textwidth]{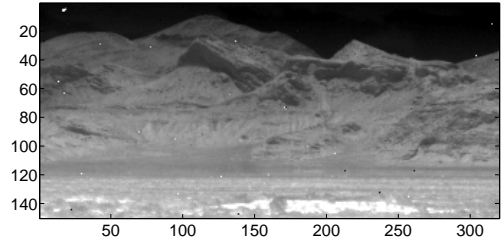}}
\end{minipage}
\caption{Contrast--enhanced spectral--mean image of the hypercube dataset.}
\label{fig:Ex3Sig}
\end{figure}

As a real life signal we consider the case of a hyperspectral image $F\in\R^{h\,\times\,v\,\times\,d}$, with $d$ frequency channels, $h\,\times\,v$ pixels $p_{ij}$ and their corresponding signatures $s_{ij}\in\R^d$. In Figure \ref{fig:Ex3Sig} we show the hyperspectral image of an area where some chemical has been released in the atmosphere. The problem is, given the spectral signature of a chemical $s_c\in\R^d$, to classify the pixels of the hypercube in order to identify the chemical position in the air \cite{manolakis2002detection}. This can be done using a classifier like the Adaptive Cosine Estimator (ACE) which assign to each pixel a value as follows:
\begin{equation}
y(p_{ij})_{ACE} = \frac{[s_{ij}^T \Sigma^{-1} s_c]^2  }{ s_c^T \Sigma^{-1} s_c s_{ij}^T \Sigma^{-1}s_{ij}}.
\label{eq:ACE}
\end{equation}
 where $\Sigma\in\R^{d\ \times\ d}$ is the covariance matrix obtained considering each pixel as an observation and each frequency as a variable.

\begin{figure}
\centering
\begin{minipage}[b]{.48\linewidth}
  \centering
  \centerline{\includegraphics[width=\linewidth]{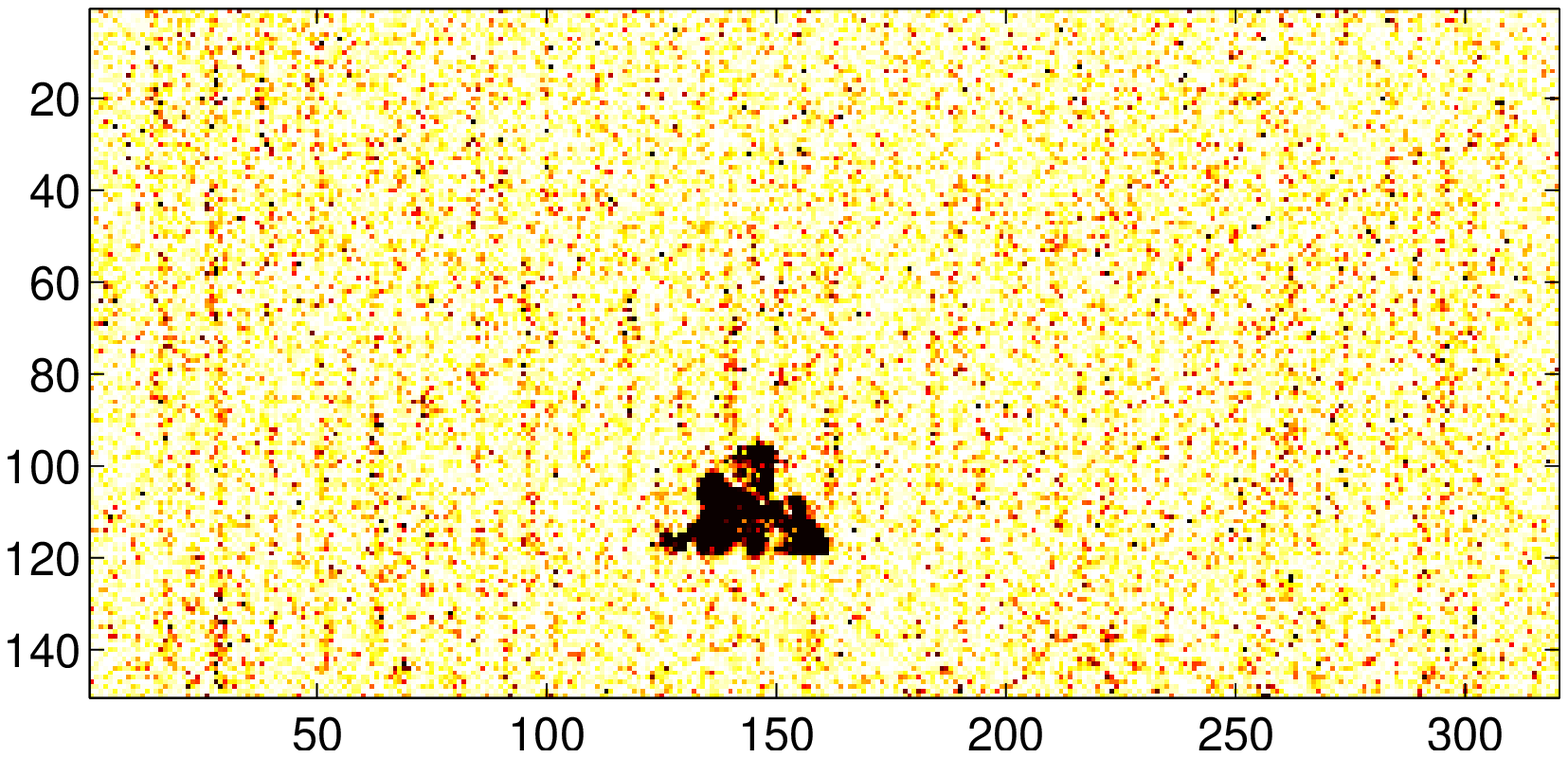}}
\end{minipage}
\hfill
\begin{minipage}[b]{0.48\linewidth}
  \centering
  \centerline{\includegraphics[width=\linewidth]{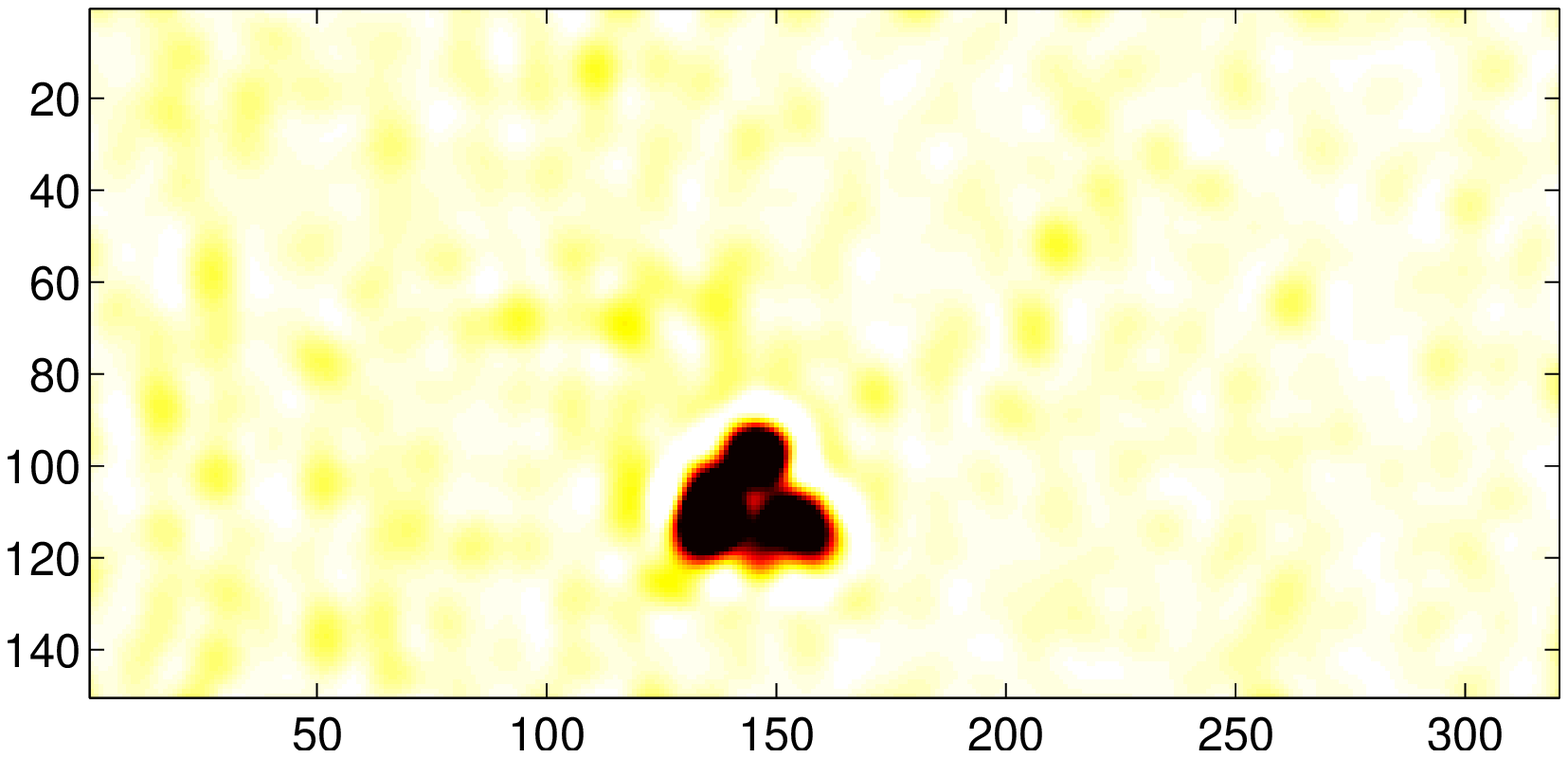}}
\end{minipage}
%
\caption{Pixel classification using the raw (left) and pre-- and post--processed (right) dataset.}\label{fig:Ex3ACE}
\end{figure}

If we apply \eqref{eq:ACE} to the raw hypercube data we obtain the classification showed in Figure~\ref{fig:Ex3ACE} (left), where the darker the pixel is the higher is the probability that it contains the chemical, while the lighter the lower is the probability. We observe that many isolated pixels are classified in black, those are all false alarms that can be due to various factors like noise or sensor malfunction.

We can use the MIF algorithm to first pre--process the hypercube, removing from each frequency channel the first IMF, and then we reapply MIF to post--process the ACE classification values, removing the first IMF produced in the decomposition. In doing so the performance of the classifier are improved as showed in Figure~\ref{fig:Ex3ACE} (right). For more details on this topic we refer the interested reader to \cite{cicone2015Hyperspectral}.

\section{Outlook}
\label{sec:outlook}

Inspired by the convergence and stability properties of Iterative Filtering and the extreme smoothness of the Fokker--Planck (FP) filters \cite{cicone2014adaptive}, in this paper we introduce the Multidimensional Iterative Filtering (MIF) algorithm for the decomposition of multidimensional non--stationary signals. Furthermore we extend FP filters to higher dimensions and we test the MIF technique, equipped with these filters, on both artificial and real life datasets. These tests show the ability of MIF to properly decompose a non--stationary signal into IMFs even of completely different nature.

This is just a preliminary work in the analysis of the MIF method. There are several open problems that need to be addressed like, for instance: Is there a simple and effective way to identify and use higher dimensional position of minima and maxima of a multidimensional signal in the derivation of the filter support? What about the convergence of MIF? Are there sufficient conditions on the higher dimensional filters that ensure the convergence of this method? What are the connections with previously developed filtering techniques like the iterative steering kernel regression method \cite{takeda2007kernel}?

Finally we point out that once a multidimensional signal is decomposed into IMFs the next step would be to apply a time--frequency analysis on each of them. Currently the most common way to perform a time--frequency analysis on 1D IMFs is via the estimation of instantaneous frequencies by means of the so called Hilbert Transform \cite{boashash1992estimating,cohen1995time,huang1998empirical}. However this approach cannot be directly  extended to handle higher dimensional IMFs, so it remains an open problem how to compute instantaneous frequency directly in higher dimensions.




\bibliography{refs4}

\begin{thebibliography}{10}

\bibitem{boashash1992estimating}
B.~Boashash.
\newblock Estimating and interpreting the instantaneous frequency of a signal.
  i. fundamentals.
\newblock {\em Proceedings of the IEEE}, 80(4):520--538, 1992.

\bibitem{cicone2014adaptive}
A.~Cicone, J.~Liu, and H.~Zhou.
\newblock Adaptive local iterative filtering for signal decomposition and
  instantaneous frequency analysis.
\newblock {\em preprint arXiv:1411.6051}, 2014.

\bibitem{cicone2015Hyperspectral}
A.~Cicone, J.~Liu, and H.~Zhou.
\newblock Hyperspectral chemical plume detection algorithms based on iterative
  filters decomposition.
\newblock {\em preprint}, 2015.

\bibitem{cohen1995time}
L.~Cohen.
\newblock {\em Time-frequency analysis}, volume 1406.
\newblock Prentice Hall PTR Englewood Cliffs, NJ:, 1995.

\bibitem{daubechies2011synchrosqueezed}
I.~Daubechies, J.~Lu, and H.-T. Wu.
\newblock Synchrosqueezed wavelet transforms: an empirical mode
  decomposition-like tool.
\newblock {\em Applied and computational harmonic analysis}, 30(2):243--261,
  2011.

\bibitem{dragomiretskiy2014variational}
K.~Dragomiretskiy and D.~Zosso.
\newblock Variational mode decomposition.
\newblock {\em IEEE transactions on signal processing}, 62(1-4):531--544, 2014.

\bibitem{gilles2013empirical}
J.~Gilles.
\newblock Empirical wavelet transform.
\newblock {\em Signal Processing, IEEE Transactions on}, 61(16):3999--4010,
  2013.

\bibitem{hou2011adaptive}
T.~Y. Hou and Z.~Shi.
\newblock Adaptive data analysis via sparse time-frequency representation.
\newblock {\em Advances in Adaptive Data Analysis}, 3(01n02):1--28, 2011.

\bibitem{huang1998empirical}
N.~E. Huang, Z.~Shen, S.~R. Long, M.~C. Wu, H.~H. Shih, Q.~Zheng, N.-C. Yen,
  C.~C. Tung, and H.~H. Liu.
\newblock The empirical mode decomposition and the hilbert spectrum for
  nonlinear and non-stationary time series analysis.
\newblock {\em Proc. of the Royal Soc. of London. Ser. A}, 454(1971):903--995,
  1998.

\bibitem{lin2009iterative}
L.~Lin, Y.~Wang, and H.~Zhou.
\newblock Iterative filtering as an alternative algorithm for empirical mode
  decomposition.
\newblock {\em Advances in Adaptive Data Analysis}, 1(04):543--560, 2009.

\bibitem{manolakis2002detection}
D.~Manolakis and G.~Shaw.
\newblock Detection algorithms for hyperspectral imaging applications.
\newblock {\em Signal Processing Mag., IEEE}, 19(1):29--43, 2002.

\bibitem{pustelnik2012multicomponent}
N.~Pustelnik, P.~Borgnat, and P.~Flandrin.
\newblock A multicomponent proximal algorithm for empirical mode decomposition.
\newblock In {\em Signal Processing Conference (EUSIPCO), 2012 Proceedings of
  the 20th European}, pages 1880--1884. IEEE, 2012.

\bibitem{selesnick2011resonance}
I.~W. Selesnick.
\newblock Resonance-based signal decomposition: A new sparsity-enabled signal
  analysis method.
\newblock {\em Sig. Proc.}, 91(12):2793--2809, 2011.

\bibitem{takeda2007kernel}
H.~Takeda, S.~Farsiu, and P.~Milanfar.
\newblock Kernel regression for image processing and reconstruction.
\newblock {\em Image Processing, IEEE Transactions on}, 16(2):349--366, 2007.

\bibitem{wei1998instantaneous}
D.~Wei and A.C. Bovik.
\newblock On the instantaneous frequencies of multicomponent am-fm signals.
\newblock {\em Signal Processing Letters, IEEE}, 5(4):84--86, 1998.

\bibitem{wu2011one}
H.-T. Wu, P.~Flandrin, and I.~Daubechies.
\newblock One or two frequencies? the synchrosqueezing answers.
\newblock {\em Adv. in Adap. Data An.}, 3(01n02):29--39, 2011.

\bibitem{wu2009ensemble}
Z.~Wu and N.~E. Huang.
\newblock Ensemble empirical mode decomposition: a noise-assisted data analysis
  method.
\newblock {\em Advances in adaptive data analysis}, 1(01):1--41, 2009.

\bibitem{wu2009multi}
Z.~Wu, N.~E. Huang, and X.~Chen.
\newblock The multi--dimensional ensemble empirical mode decomposition method.
\newblock {\em Adv. in Adap. Data An.}, 1(03):339--372, 2009.

\bibitem{yang2014synchrosqueezed}
H.~Yang and L.~Ying.
\newblock Synchrosqueezed curvelet transform for two--dimensional mode
  decomposition.
\newblock {\em SIAM Journal on Mathematical Analysis}, 46(3):2052--2083, 2014.

\end{thebibliography}
\bibliographystyle{plain}

\end{document}